\documentclass[paper=letter, fontsize=12pt]{article}
\usepackage[english]{babel} % English language/hyphenation
\usepackage{amsmath,amsfonts,amsthm} % Math packages
\usepackage[utf8]{inputenc}
\usepackage{float}
\usepackage{lipsum} % Package to generate dummy text throughout this template
\usepackage{blindtext}
\usepackage{graphicx} 
\usepackage{caption}
\usepackage[sc]{mathpazo} % Use the Palatino font
\usepackage[T1]{fontenc} % Use 8-bit encoding that has 256 glyphs
\usepackage{microtype} % Slightly tweak font spacing for aesthetics
\usepackage[hmarginratio=1:1,top=32mm,columnsep=20pt]{geometry} % Document margins
\usepackage{multicol} % Used for the two-column layout of the document
\usepackage{booktabs} % Horizontal rules in tables
\usepackage{float} % Required for tables and figures in the multi-column environment - they need to be placed in specific locations with the [H] (e.g. \begin{table}[H])
\usepackage{hyperref} % For hyperlinks in the PDF
\usepackage{lettrine} % The lettrine is the first enlarged letter at the beginning of the text
\usepackage{paralist} % Used for the compactitem environment which makes bullet points with less space between them
\usepackage{abstract} % Allows abstract customization
\usepackage{titlesec} % Allows customization of titles
\usepackage{natbib, url}
\usepackage{subfig}
\usepackage{listings}
\renewcommand\thesection{\Roman{section}} % Roman numerals for the sections
\renewcommand\thesubsection{\Roman{subsection}} % Roman numerals for subsections

\titleformat{\section}[block]{\large\scshape\centering}{\thesection.}{1em}{} % Change the look of the section titles
\titleformat{\subsection}[block]{\large}{\thesubsection.}{1em}{} % Change the look of the section titles
\usepackage{fancyhdr} % Headers and footers
\graphicspath{{./pic/}}

\title{\vspace{-15mm}\fontsize{24pt}{10pt}\selectfont\textbf{A Note on the comparison of Nearest Neighbor Gaussian Process (NNGP) based models}} % Article title

\author{
\large
{\textsc{Lu Zhang}}\\[2mm]
{\textsc{UCLA Department of Biostatistics }}\\[2mm]
\normalsize \href{mailto:lu.zhang@ucla.edu}{Lu.Zhang@ucla.edu}\\[2mm] 
\large
{\textsc{Sudipto Banerjee}}\\[2mm]
{\textsc{UCLA Department of Biostatistics }}\\[2mm]
\normalsize \href{mailto:sudipto@ucla.edu}{sudipto@ucla.edu}\\[2mm] 
}
\date{Nov 08, 2018}
\providecommand{\keywords}[1]{\textbf{\textit{Key words:}} #1}
\begin{document}
\maketitle % Insert title
\thispagestyle{fancy} % All pages have headers and footers

\label{firstpage}

\begin{abstract}
	  This note is devoted to the comparison between two Nearest-neighbor Gaussian processes (NNGP) based models: the response NNGP model and the latent NNGP model. We exhibit that the comparison based on the Kullback-Leibler divergence (KL-D) from the NNGP based models to their parent GP based model can result in reverse conclusions in different parameter spaces. And we suggest a heuristic explanation on the phenomenon that the latent NNGP model tends to outperform the response NNGP model in approximating their parent GP based model.
\end{abstract}
\newpage
\section{Introduction}\label{sec: intro}
Big data problems in Geostatistics have generated substantial interest over the last decade. Gaussian Processes (GPs) based models are among the most popular for modeling location-referenced spatial data. Fitting GP based models to large datasets, however, requires matrix operations that quickly become prohibitive (see, e.g., \citep{banerjee2017high}. An early approach of approximating the likelihood of GP based models \citep{ve88} has attracted much attention recently due to its proven effectiveness (also see \citep{stein2014}. More recently, innovated by Vecchia's approach, \citep{datta16} proposed an approach of constructing a GP called Nearest Neighbor Gaussian Process (NNGP) based upon a given GP (parent GP). The merit of NNGP is that it introduces sparsity in the precision matrix of its finite realization, which lessens the storage requirement and computational burden in the model fitting, while the inference obtained from an NNGP based model is comparable to its parent GP based model. Some of the practical implementations of the NNGP is discussed in \cite{zhang2018practical}, while \citep{katzfuss2017general} discuss certain generalizations. 

In this note, we focus on discussing the performance of two NNGP based models: the response NNGP and latent NNGP models. %We use Kullback-Leibler divergence (KL-D) as the criterion in the performance comparison. 
The note is organized as follows. First, we introduce the Kullback-Leibler divergence (KL-D) and discuss the KL-D for Gaussian models with latent variables. Second, we discuss a claim from \citep{katzfuss2017general} which compares the response NNGP and latent NNGP models based on KL-D. Third, we suggest a heuristic explanation on the phenomenon that the latent NNGP model tends to outperform the response NNGP model in approximating their parent GP based model

\section{KL-D of Gaussian models with latent variables}\label{sec: KL-D intro}

\noindent The Kullback-Leibler divergence (KL-D) is widely used as a measure of the difference between two probability distributions. The KL-D from probability distribution $Q$ to probability distribution $P$ is defined as

\begin{equation}
D_{KL}(P || Q) = \int \log{\frac{dP}{dQ}}dP\; , \;
\end{equation}
where $\frac{dP}{dQ}$ is the Radon-Nikodym derivative of $P$ with respect to $Q$. In particular, for two Gaussian distributions $P, Q$ on $\mathcal{R}^N$ with mean vector $\mu_P, \mu_Q$ and covariance matrix $\Sigma_P, \Sigma_Q$, respectively, the KL-D from $Q$ to $P$ is given by:

\begin{equation}\label{eq: KL-D_Gaussian}
D_{KL}(P || Q) = -\frac{1}{2}\{\log{det\Sigma_P} - \log{det\Sigma_Q} + N - tr(\Sigma_Q^{-1}\Sigma_P) - (\mu_P - \mu_Q)^{\top} \Sigma_Q^{-1} (\mu_P - \mu_Q)\}
\end{equation}

A KL-D of zero indicates that $P, Q$ satisfy $P = Q$ almost everywhere. The higher the KL-D is, the more differently the two distributions behave. Thus, the value of KL-D is popularly recognized as a criterion in comparing the performance of models. A valid comparison based on KL-D requires that the distributions based on models are all defined on the same space. For hierarchical models, latent variables are often introduced to model unobserved variables, which results in a distribution on the joint space of the latent and observed variables. It is common to use the joint distribution of the latent and observed variables to calculate KL-D for theoretical comparison. However, in practice, the hierarchical models with latent variables are not preferred for obtaining inference in Markov Chain Monte Carlo (MCMC) algorithms. The latent variables are treated as parameters to be sampled since they are unobserved, but needed in each iteration of MCMC algorithms. This enlarges the dimension of the parameter space and dramatically lowers the sampling efficacy in MCMC algorithms. It is popular to integrate out latent variables and use the marginal distribution on the collapsed space to obtain the inference in MCMC algorithms. While one model may be better than another based on the KL-Ds on the joint space of observed and latent variables, the question remains whether the model will still be better on the collapsed space of observed variables. The following example shows that marginalization neither preserves volume nor keeps the order of KL-D.

Assume the true model is $y | w \sim N(w, 1)$ with latent variable $w \sim N(0, 1)$, now consider two models 

\begin{equation}
\begin{aligned}
model1: y | w &\sim N(w, 0.5)\; , \; w \sim N(0, 0.5) \; \\  
model2: y | w &\sim N(w, 2.5)\; , \; w \sim N(0, 0.5)
\end{aligned}
\end{equation} 
The true distribution $P$ and the probability distributions for two models $Q_1, Q_2$ follows

\begin{equation}
\begin{aligned}
P(y, w) \sim N \left(\begin{bmatrix}0\\0\end{bmatrix}, \begin{bmatrix}2.0 & 1.0\\ 1.0 & 1.0 \end{bmatrix} \right) ; \;  &P(y) \sim N(0, 2.0)\\
Q_1(y, w) \sim N \left(\begin{bmatrix}0\\0\end{bmatrix}, \begin{bmatrix}1.0 & 0.5\\ 0.5 & 0.5 \end{bmatrix} \right)  ; \; &Q_1(y) \sim N(0, 1.0)\\
Q_2(y, w) \sim N \left(\begin{bmatrix}0\\0\end{bmatrix}, \begin{bmatrix}3.0 & 0.5\\ 0.5 & 0.5 \end{bmatrix} \right)  ; \; &Q_2(y) \sim N(0, 3.0)
\end{aligned}
\end{equation}
Then through \ref{eq: KL-D_Gaussian}, the KL-D from $Q1$ to $P$ and the KL-D from $Q2$ to $P$ satisfy:

\begin{equation}\label{eq: KL-D_com1}
D_{KL}(P(y, w) || Q_1(y, w)) < D_{KL}(P(y, w) || Q_2(y, w)) ; \; D_{KL}(P(y) || Q_1(y)) > D_{KL}(P(y) || Q_2(y))
\end{equation}

This indicates that $Q_1$ is closer to $P$ than $Q_2$ on the joint space of observed and latent variables, while $Q_1$ is further from $P$ than $Q_2$ on the collapsed space of observed variables. Now, let the second model be: 

\begin{equation}
model2: y | w \sim N(w, 1.5)\; , \; w \sim N(0, 1.5)
\end{equation}
With the correponding $Q_2$

\begin{equation}
Q_2(y, w) \sim N \left(\begin{bmatrix}0\\0\end{bmatrix}, \begin{bmatrix}3.0 & 1.5\\ 1.5 & 1.5 \end{bmatrix} \right) ; \; Q_2(y) \sim N(0, 3.0)
\end{equation}
Then the KL-D comparison becomes:

\begin{equation}
D_{KL}(P(y, w) || Q_1(y, w)) > D_{KL}(P(y, w) || Q_2(y, w)); \; D_{KL}(P(y) || Q_1(y)) > D_{KL}(P(y) || Q_2(y))\; ,
\end{equation}
which shows that $Q_1$ is closer to $P$ than $Q_2$ on both the joint and the collapsed spaces. The above example illustrates that the comparison made on the augmented space may not hold on the collapsed space when using KL-Ds. In the next section, we will show that a claim in \citep{katzfuss2017general} may fail when we use KL-Ds on a collapsed space.

\section{The comparison of the performance of NNGP based models} \label{Sec: compa_KL-D}
The preceding section illustrates the comparison of KL-Ds with a trivial hierarchical model. In this section, we discuss the performance of NNGP based models using the KL-Ds from them to their parent GP based models. Consider $\{y(s): s \in \mathcal{D}\}$ be the process of interest over domain $\mathcal{D} \subset R^d, d \in N^+$, and $y(s)$ can be decomposed as $y(s)  = w(s) + \epsilon(s)$ for some latent process $w(s)$ and white noise process $\epsilon(s)$. Assume that $y(s)$ and $w(s)$ are known to follow certain GPs. We use an NNGP derived from the true GP in modeling as an alternative of the true GP. Depending on which process NNGP is assigned for, there will be a response NNGP model and a latent NNGP model. The former assigns NNGP for $y(s)$ and the latter assigns NNGP for $w(s)$. Using a Directed Acyclic Graph (DAG) built on an augmented latent space, \citep{katzfuss2017general} show that the KL-D from the latent NNGP model to the true model is no more than that of the response NNGP model, suggesting that the latent NNGP model is better than the response NNGP model. However, as pointed out in the last section, the claim based on KL-Ds on an augmented space is not guaranteed to hold on a collapsed space. Here we provide an example with numerical results to show a response NNGP model that might outperform a latent NNGP model on a collapsed space.  

Assume the observed location set is $S = \{s_1, s_2, s_3\}$, $w(s_1, s_2, s_3)$ has covariance matrix $\sigma^2 R$ with correlation matrix:

\begin{equation}
R = \begin{bmatrix}
1& \rho_{12}  & \rho_{13}\\ 
\rho_{12}& 1 & \rho_{23} \\ 
\rho_{13}& \rho_{23} & 1
\end{bmatrix}
\end{equation}

Let us suppress the connection between knots $s_2$ and $s_3$ in the DAG corresponding to the finite realization of the NNGP on $S$. Then the covariance matrix of $y(s_1, s_2, s_3)$ of the response NNGP model $\Sigma_R$ and that of the latent NNGP model $\Sigma_l$ have the following forms:

\begin{equation}
\Sigma_R = \sigma^2 \begin{bmatrix}
1 + \delta^2 & \rho_{12}  & \frac{\rho_{12} \rho_{23}}{ 1 + \delta^2}\\ 
\rho_{12}& 1 + \delta^2 & \rho_{23} \\ 
\frac{\rho_{12} \rho_{23}}{ 1 + \delta^2}& \rho_{23} & 1 + \delta^2
\end{bmatrix}\; , \;
\Sigma_l = \sigma^2 \begin{bmatrix}
1 + \delta^2 & \rho_{12}  & \rho_{12} \rho_{23}\\ 
\rho_{12}& 1 + \delta^2 & \rho_{23} \\ 
\rho_{12} \rho_{23}& \rho_{23} & 1 + \delta^2 \; ,
\end{bmatrix}
\end{equation}

where $\delta^2 = \frac{\tau^2}{\sigma^2}$ is the noise-to-signal ratio with $\tau^2$ as the variance of the noise process $\epsilon(s)$. By the sufficient and necessary condition of $R$ being positive-definite:

\begin{equation}
1- (\rho_{12}^2 + \rho_{13}^2 + \rho_{23}^2) + 2\rho_{12}\rho_{13}\rho_{23} > 0 \; ,\; 1 - \rho_{12}^2 > 0\; ,
\end{equation}
it is easy to show that $\Sigma_R$ and $\Sigma_l$ are also positive-definite. If $\rho_{13} = \frac{\rho_{12} \rho_{23}}{ 1 + \delta^2}$, then the KL-D from the response NNGP model to the true model always equals zero, which is no more than the KL-D from the latent NNGP model to the true model. If $\rho_{13} = \rho_{12} \rho_{23}$, then the KL-D of the latent NNGP model to the true model always equals zero, which reverses the relationship. More examples were found in the trial and error of this study. The numerical studies can be found in \url{https://luzhangstat.github.io/notes/KL-D_com.html}

While there can be found examples for both sides as discussed above, we observed that the latent NNGP model is in general a better approximation to the parent GP than the response NNGP model based on the KL-Ds when fixing process parameters at the true value. Though the conditions for each sides is still unclear for the author, we try to give a heuristic explanation in next section to show why latent NNGP model tends to be better than the response NNGP. 

\section{A discussion on the possible reason behind the observed relationship}

Let the covariance matrix of $y(S)$ of the parent GP based models be $C + \tau^2I$ where $C$ is the covariance matrix of the latent process $w(S)$. Define the Vecchia approximation of the precision matrix $C^{-1}$ and $K^{-1} = \{C + \tau^2 I\}^{-1}$ be
\begin{equation}
Vecchia(C^{-1}) = {\tilde{C}}^{-1}\; , \; Vecchia(K^{-1}) = \tilde{K}^{-1}
\end{equation}
so that the covariance matrix of $y(S)$ of the latent NNGP model is $\tilde{C} + \tau^2 I$ and the precision matrix of $y(S)$ of the response NNGP model is $\tilde{K}^{-1}$. We denote the error matrix of the Vecchia approximation of $C^{-1}$ be $E$. $E$ is small so that $\tilde{C}^{-1}$ approximates $C^{-1}$ well. With the same observed location $S$ and the fixed number of nearest neighbors, the error matrix of the Vecchia approximation of $K^{-1}$ is believed to be close to $E$
\begin{equation}
C^{-1} = \tilde{C}^{-1} + E\; ; \; K^{-1} = \tilde{K}^{-1} + \mathcal{O}(E).
\end{equation}
Represent the precision matrix of $y(S)$ of the parent GP based model and the latent NNGP model to be
\begin{equation}
\begin{aligned}
(C+\tau^2I)^{-1} &= C^{-1} - C^{-1}M^{-1}C^{-1}\; , M = C^{-1} + \tau^{-2} I\\
(\tilde{C} + \tau^2I)^{-1} &= \tilde{C}^{-1} - \tilde{C}^{-1} M^{\ast-1}\tilde{C}^{-1}\; ,
M^\ast = \tilde{C}^{-1} + \tau^{-2}I
\end{aligned}
\end{equation}
Check the difference 
\begin{align*}
(C+\tau^2I)^{-1} - (\tilde{C}+ \tau^2I)^{-1} 
&= C^{-1} - C^{-1}M^{-1}C^{-1} - \tilde{C}^{-1} + \tilde{C}^{-1} M^{\ast-1}\tilde{C}^{-1}\\
&= \underbrace{E - EM^{-1}\tilde{C}^{-1} - \tilde{C}^{-1}M^{-1}E - \tilde{C}^{-1}(M^{-1} - M^{\ast-1})\tilde{C}^{-1}}_{B} - \underbrace{EM^{-1}E}_{\mathcal{O}(E^2)}
\end{align*}
%If the leading term $B \approx \mathcal{O}(E^{1 + \delta})$ for some $\delta > 0$ then latent NNGP model will outperform the response NNGP in general cases. This strong condition is still to be explored. However, 
We find that the leading term $B$ tends to shink the error matrix from the Vecchia's approximation of $C^{-1}$. Consider representing $B$ in terms of $\tilde{C}^{-1}$, $M^\ast$ and $E$, where $E$ is assumed to be nonsingluar. 

%\begin{align*}
%M^{-1} &= (M^\ast + E)^{-1} = M^{\ast-1} - M^{\ast-1}(E^{-1} + M^{\ast-1})^{-1}M^{\ast-1}\\
%M^{-1} - M^{\ast -1} &= - M^{\ast-1}(E^{-1} + M^{\ast-1})^{-1}M^{\ast-1}\\
%EM^{-1}\tilde{C^{-1}} &= E[M^{\ast-1} - M^{\ast-1}(E^{-1} + M^{\ast-1})^{-1}M^{\ast-1}]\tilde{C^{-1}}\\
%&= EM^{\ast-1}\tilde{C^{-1}} - EM^{\ast-1}(E^{-1} + M^{\ast-1})^{-1}M^{\ast-1}\tilde{C^{-1}}\\
%\tilde{C^{-1}}M^{-1}E &= \tilde{C^{-1}}M^{\ast-1}E - \tilde{C^{-1}}M^{\ast-1}(E^{-1} + M^{\ast-1})^{-1}M^{\ast-1}E
%\end{align*}

\begin{equation}
\label{eq: B_1}
\begin{aligned}
B =& \;E - EM^{\ast-1}\tilde{C}^{-1} + EM^{\ast-1}(E^{-1} + M^{\ast-1})^{-1}M^{\ast-1}\tilde{C}^{-1} \\
&\;- \tilde{C}^{-1}M^{\ast-1}E + \tilde{C}^{-1}M^{\ast-1}(E^{-1} + M^{\ast-1})^{-1}M^{\ast-1}E
+\tilde{C}^{-1}M^{\ast-1}(E^{-1} + M^{\ast-1})^{-1}M^{\ast-1}\tilde{C}^{-1}
\end{aligned}
\end{equation}
By Sherman–Morrison–Woodbury formula and the expansion $(I + X)^{-1} = \sum_{n = 0}^{\infty} \{-X\}^{n}$, we have
\begin{align*}
(E^{-1} + M^{\ast-1})^{-1}M^{\ast-1} &= \{M^{\ast}(E^{-1} + M^{\ast-1})\}^{-1} = \{M^\ast E^{-1} + I\}^{-1} \\
& = I - \{I + E M^{\ast-1}\}^{-1} = I - \{I - EM^{\ast-1} + \mathcal{O}(E^2)\} \\
&= EM^{\ast-1} + \mathcal{O}(E^2)
\end{align*}
Use the above equations and exclude the term in order $\mathcal{O}(E^2)$ in the expression of $B$, we have the leading term of the difference to be
\begin{equation}
%E + \tilde{C}^{-1}M^{\ast-1}EM^{\ast-1}\tilde{C}^{-1} - EM^{\ast-1}\tilde{C}^{-1} - \tilde{C}^{-1}M^{\ast-1}E
B = (I - \tilde{C}^{-1}M^{\ast-1}) E (I - M^{\ast-1}\tilde{C}^{-1})  =
(I + \tau^2 \tilde{C}^{-1})^{-1}E(I + \tau^2 \tilde{C}^{-1})^{-1} 
\end{equation}
By the spectrum decomposition theorem, there exist an orthogonal matrix $P$ such that 
\begin{align*}
(I + \tau^2 \tilde{C}^{-1}) = P^\top (I + \tau^2 D)P
\end{align*}
where $D$ is a diagonal matrix whose elements on diagonal are positive. Thus
\begin{equation}
\begin{aligned}
||B||_F &= ||P^\top (I + \tau^2 D)^{-1}PEP^\top(I + \tau^2 D)^{-1}P||_F = || (I + \tau^2 D)^{-1}PEP^\top(I + \tau^2 D)^{-1}||_F\\
 &\leq || PEP^\top||_F = ||E||_F
\end{aligned}
\end{equation}
where $||\cdot||_F$ is the Frobenius norm. The inequality also holds for the absolute value of determinant and $p$ norms. And the equality holds if and only if $\tau^2 = 0$ when the difference is the same as the error matrix for response NNGP model. Thus we conclude that the latent model shrinks the error from the Vecchia approximation. And this might be a reason why the latent NNGP model in general outperforms the response NNGP model in approximating the parent GP based model as observed in the studies based on KL-Ds.

\appendix

\bibliographystyle{ba}  
\bibliography{lubib} 

\begin{thebibliography}{6}
\newcommand{\enquote}[1]{``#1''}
\expandafter\ifx\csname natexlab\endcsname\relax\def\natexlab#1{#1}\fi
\expandafter\ifx\csname url\endcsname\relax
  \def\url#1{{\tt #1}}\fi
\expandafter\ifx\csname urlprefix\endcsname\relax\def\urlprefix{URL }\fi
\ifx\endbibitem\undefined \let\endbibitem\relax\fi

\bibitem[{Banerjee(2017)}]{banerjee2017high}
Banerjee, S. (2017).
\newblock \enquote{High-Dimensional Bayesian Geostatistics.}
\newblock {\em Bayesian Analysis\/}, 12: 583--614.
\endbibitem

\bibitem[{Datta et~al.(2016)Datta, Banerjee, Finley, and Gelfand}]{datta16}
Datta, A., Banerjee, S., Finley, A.~O., and Gelfand, A.~E. (2016).
\newblock \enquote{Hierarchical Nearest-Neighbor Gaussian Process Models for
  Large Geostatistical Datasets.}
\newblock {\em Journal of the American Statistical Association\/}, 111:
  800--812.
\newline\urlprefix\url{http://dx.doi.org/10.1080/01621459.2015.1044091}
\endbibitem

\bibitem[{Katzfuss and Guinness(2017)}]{katzfuss2017general}
Katzfuss, M. and Guinness, J. (2017).
\newblock \enquote{A General Framework for Vecchia Approximations of Gaussian
  Processes.}
\newblock {\em arXiv preprint arXiv:1708.06302\/}.
\endbibitem

\bibitem[{Stein(2014)}]{stein2014}
Stein, M.~L. (2014).
\newblock \enquote{Limitations on Low Rank Approximations for Covariance
  Matrices of Spatial Data.}
\newblock {\em Spatial Statistics\/}, 8(0): 1--19.
\endbibitem

\bibitem[{Vecchia(1988)}]{ve88}
Vecchia, A.~V. (1988).
\newblock \enquote{Estimation and Model Identification for Continuous Spatial
  Processes.}
\newblock {\em Journal of the Royal Statistical society, Series B\/}, 50:
  297--312.
\endbibitem

\bibitem[{Zhang et~al.(2018)Zhang, Datta, and Banerjee}]{zhang2018practical}
Zhang, L., Datta, A., and Banerjee, S. (2018).
\newblock \enquote{Practical Bayesian Modeling and Inference for Massive
  Spatial Datasets On Modest Computing Environments.}
\newblock {\em arXiv preprint arXiv:1802.00495\/}.
\endbibitem

\end{thebibliography}

\label{lastpage}	

\end{document}